\documentclass[a4paper,11pt]{amsart}
 \usepackage{amsfonts,mathrsfs}
 \usepackage{amsthm}
 \usepackage{amssymb}
 \usepackage{enumerate}
 \usepackage{tikz-cd}
 \usepackage{MnSymbol}
 \usepackage{csquotes}
 \usepackage{cleveref}
 \usepackage{paralist}
 \usepackage[]{xfrac}
 \usepackage[]{url}

\theoremstyle{definition}

\theoremstyle{plain}

\newtheorem*{thm*}{Theorem}

\newtheorem*{cor*}{Corollary}

\newtheorem*{con*}{Conjecture}

\title[]{Landau on Chess Tournaments\\ and {Google}'s {PageRank}}
\author{Rainer Sinn}
\address{Rainer Sinn: Universität Leipzig, Mathematisches Institut, Augustusplatz 10, 04109 Leipzig, Germany} 
\email{rainer.sinn@uni-leipzig.de}
\author{Günter M.~Ziegler}
\address{Günter M.~Ziegler: Freie Universität Berlin, Institut für Mathematik, Arnimallee 2, 14195 Berlin, Germany}
\email{guenter.m.ziegler@fu-berlin.de}

\begin{document}


\begin{abstract}
In his first mathematical paper, published in 1895 when he was 18, Edmund Landau suggested a new way to determine the winner of a chess tournament by not simply adding for each player the fixed number of points they would get for each win or draw, but rather by considering the performance of all players in the tournament relative to each other: each player would get more credit for games won or drawn against stronger players. Landau called this \enquote{relative Wertbemessung}, which translates to relative score. The basic idea from linear algebra behind this scoring system was rediscovered and reused in many contexts since 1895; in particular, it is a central ingredient in \emph{Google}'s \emph{PageRank}.
\end{abstract}

\maketitle

\section*{Landau and Chess}
At the end of the 19th century, there was a lively discussion in the world of chess about how to determine the winner of a round robin chess tournament. 
So the point of departure is a tournament where each player plays once against each other player and the possible outcomes are a win for one side (and thus a loss for the opponent), or a draw. 
In the common scoring system for chess tournaments at the time,
every player would get $1$ point for each win (in which case the opponent gets $0$ points), and in case of a draw both players get $\sfrac{1}{2}$ point.
Then one would add for each player all his points --- and the player with the highest sum wins. Of course several players might have the same maximal score, so the winner might not be unique. 
Such systems are still widely used, also in other sports like American football and European football (soccer), where in each season each team plays each other team twice. In European football, each team gets 3 points for a win and 1 point for a draw.

The question debated in chess journals at the time was whether this overall score accurately reflects the elusive \emph{quality} of the players at the tournament. Should a player who performed very well against \enquote{strong} players 
(and less well against \enquote{weaker} opponents) be considered stronger than someone who performed badly against the strong players but managed to scrape more \enquote{easy} points against \enquote{weaker} players?

In his first ever publication, \enquote{Zur relativen Wertbemessung der Turnierresultate} in the magazine
\textit{Deutsches Wochenschach} \cite{Wochenschach}, which was edited by Siegbert Tarrasch, 
a world class chess player, 
the young Edmund Landau suggested a new way to compute an overall score from the game results that should better capture the quality of the performance in a tournament. Thus he brought a mathematical perspective into an on-going discussion, which he expanded in a follow-up paper in 1914~\cite{Preisverteilung}.

\section*{Relative Scores and Eigenvalues}
The point of departure for Landau's method is the \emph{game results matrix}. For a round robin tournament with $n$ players $1,\ldots,n$, this is an $(n\times n)$-matrix $A$ whose entry in row $i$ and column $j$ records the points that player~$i$ scores in the game against player~$j$. Thus the $i$th row records the points won by player~$i$. For example, the matrix
\[
	A \ =\ 
  \begin{pmatrix}
    0 & 1 & \sfrac12 \\
    0 & 0 & 1 \\
    \sfrac12 & 0 & 0 
  \end{pmatrix}
\]
records the results of a $3$ player round robin tournament, where player~$1$ has won against player~$2$ and played a draw against player~$3$, while player~$2$ won against player~$3$. The diagonal entries are $0$ by convention (as the players don't play against themselves); by construction $A+A^t$ also has $0$s on the diagonal, but $1$ in all off-diagonal entries.

The common way to determine the overall score of player~$i$ at the time of Landau's first article (and still common in chess today) is to compute the sum of all entries in row $i$, which gives the sum of all game results for that player. The players are then ranked according to these row sums. 

Landau's idea was to introduce the unknown \emph{quality} of each player $i$ as a variable $x_i$, assuming that the qualities determine the outcomes of the games as an expected value: Writing $a_{ij}$ for the points that player~$i$ got for his game against player~$j$, the quality $x_i$ of player~$i$ in this tournament is proportional to 
the weighted row sum
\[
  \sum_{j=1}^n a_{ij}x_j,
\]
which does not depend on $x_i$, as $a_{ii}=0$.

The goal is to now determine the unknown qualities $x_j$ from this system of equations. In matrix notation, we are solving
\[ 
	A\mathbf{x} = \lambda \mathbf{x}, 
\]
where $\mathbf{x}$ is the vector of player qualities and $\lambda$ is the proportionality factor (which is unknown as well). So in fact, this becomes an eigenvalue problem and the proposed overall score of each player is given by the entries of the unknown eigenvector $\mathbf{x}$. 

In his publication \cite{Wochenschach} in 1895, Landau first solves for $x_1,\ldots,x_n$ as functions of $\lambda$ (in fact for the quotients $x_i/x_n$ for $i=1,\ldots,n-1$ using only $n-1$ of the above linear equations) and considers the problem essentially solved at that point. He does not use the word “eigenvalue” in \cite{Wochenschach}, and he does not comment on the fact that $\lambda$ is not uniquely determined. Neither does he worry about the existence of a real eigenvalue needed for a real eigenvector whose entries can then give a ranking of the players. 

In his follow-up paper \cite{Preisverteilung} in 1914, he solves these issues with the help of a result, hot off the press at the time, which we nowadays know as the \emph{Perron--Frobenius Theorem}, and which roughly says the following: Any matrix with nonnegative real entries has a largest real eigenvalue and a corresponding eigenvector, unique up to scaling, whose entries are all nonnegative. This is the eigenvalue that Landau takes to compute the overall tournament score via a corresponding eigenvector recording all the qualities $x_i$. Taking the largest eigenvalue ensures that the score of every player is a nonnegative real number (by Perron--Frobenius).%
\footnote{To be precise, we need to assume that the matrix is irreducible. A matrix is reducible if its rows and columns can be permuted in such a way that the matrix is block upper triangular with square blocks on the diagonal. Otherwise, it is irreducible. Observe that a game results matrix can be reducible!}

Why was this system not widely adopted in competitive sports? Landau discusses his method in detail and also explains why he would not argue for this specific scoring system. His central point is that in order to maximize the score, it is not always a good idea to maximize each game result (depending on the current standing). Landau gives the following two examples of a round robin tournament with three players:%
\footnote{Note that both these matrices are reducible! (Indeed, if some player wins all his games, as player~$1$ does in these examples, then the matrix will be reducible.) Landau makes the game results matrices irreducible by deviating from actual chess game scores and considering a result $1-\varepsilon$ for player~$1$ vs.~$3$ and symmetrically a result of $\varepsilon$ for $3$ vs.~$1$ and then makes a limiting argument in terms of $\varepsilon$.}
\[
  A_1\ =\ \begin{pmatrix}
    0 & 1 & 1 \\ 
	0 & 0 & \sfrac12 \\ 
	0 & \sfrac12 & 0
  \end{pmatrix} 
  \qquad\text{and}\qquad 
  A_2\ =\ \begin{pmatrix}
    0 & 1 & 1 \\ 
	0 & 0 & 1 \\ 
	0 & 0 & 0 
  \end{pmatrix}.
\]
In both cases, player~$1$ wins both games. In the tournament $A_1$, player~$2$ draws against player~$3$, whereas he wins against player $3$ in the tournament recorded by~$A_2$. 
For the tournament $A_1$ we get that $(4,1,1)$ is an eigenvector corresponding to the dominant eigenvalue $\sfrac12$, so that players~$2$ and~$3$ have the same positive score. For $A_2$, a corresponding eigenvector is $(1,0,0)$ (all eigenvalues are $0$), so the components for players~$2$ and~$3$ are still the same, even though player~$2$ won the game against player~$3$: indeed, they both have a score of~$0$. If we take the score as their share of the prize money, players~$2$ and~$3$ would win some money in the tournament with results given by~$A_1$ (namely each $\sfrac16$ of the overall prize pool), but they both get nothing if the results are given by~$A_2$.
So winning against player $3$ here does not pay off for player $2$.
Landau concludes that such phenomena make this unsuitable for a scoring system in competitive sports.

Landau's scoring system has been applied for fun to College Football by James P.~Keener \cite{Keener} in order to popularize the Perron--Frobenius Theorem among students. In Germany, Dirk Frettl\"oh has done the same by recomputing the score of the football teams in the Bundesliga \cite{Frettloeh} for the seasons 2001/2002 until 2011/2012. In one out of the $11$ seasons in this list, the winner according to Landau's method would have been a different team. 

Landau's idea to determine the score by quality was picked up many times and indeed this is a central ingredient in \emph{Google}'s scoring of web pages by what they call~\emph{PageRank}. To apply the idea of \enquote{quality} to the internet, there are mainly two questions to address: How could or should one assign scores to websites? For this \emph{Google} would build a matrix analogous to the game results matrix. And how can one effectively compute the eigenvalue for very large matrices like the \emph{Google} matrix for the internet? Historically, the answer to the second question came first.

\section*{The Power Method}
In 1955, Maurice Kendall \cite{Kendall} considered the question of how to determine the most “popular” person in a group (or the highest-ranking element in a finite set) based on pairwise comparisons. This is similar to an overall tournament score based on chess games: The result of each game is a pairwise comparison between the players and we want to determine who is the \enquote{best} player. Kendall, however, arrived at this question through voting systems. Kendall's approach, which built on the unpublished 1952 PhD thesis 
\emph{Algebraic foundations of ranking theory} by Teh-Hsing Wei at the University of Cambridge, leads essentially to what we nowadays call the \emph{power method}. Given a game results matrix $A$, Wei suggested to compute a new score for each player by assigning to him the score of every player he has beaten and half the score of every player with whom he has drawn. Translated to matrix language, this amounts to computing $A( A \mathbf{1})$, where $\mathbf{1}$ is the vector whose entries are all equal to~$1$. 

Indeed, $A \mathbf{1}$ is the vector of scores and the next multiplication with $A$, whose entries are $1$ for a victory and $\sfrac{1}{2}$ for a draw, gives exactly the result proposed by Wei. 
Kendall takes this one step further: Why should we stop here? We can iterate this procedure, compute the sequence $\mathbf{x}_k = A^k \mathbf{1}$, and consider its limit. Again, under the mild assumptions needed for the Perron--Frobenius Theorem, the normalized sequence $y_k = \mathbf{x}_k/\|\mathbf{x}_k\|$ has such a limit. And this limit is in fact the same result as what we get from Landau's method! The power method converges to an eigenvector of the dominant eigenvalue (if the starting  vector is not orthogonal to this eigenspace). Heuristically, write our input vector $\mathbf{1} = \mathbf{x}_0$ as a linear combination of eigenvectors of $A$, say $x_0 = \sum a_i \mathbf{v}_i$ with $A\mathbf{v}_i = \lambda_i \mathbf{v}_i$. Then $\mathbf{x}_k = A^k\mathbf{x}_0 = \sum a_i \lambda_i^k \mathbf{v}_i = \lambda_1^k  \left( a_1\mathbf{v}_1 + \sum a_i \lambda_i^k / \lambda_1 \mathbf{v}_i\right)$ and the angle between this vector and $\mathbf{v}_1$ goes to $0$ (if $\lambda_1$ is the largest eigenvalue of $A$).

Kendall considers the following example illustrating how this scoring method changes the overall rating of players: his game results matrix is%
\footnote{For a game results matrix of a chess tournament, the diagonal entries would be $0$. This does not change the point of this example very much.}
\[
  A = \begin{pmatrix}
    \sfrac12 & 1 & 1 & 0 & 1 & 1 \\
    0 & \sfrac12 & 0 & 1 & 1 & 0 \\
    0 & 1 & \sfrac12 & 1 & 1 & 1 \\
    1 & 0 & 0 & \sfrac12 & 0 & 0 \\
    0 & 0 & 0 & 1 & \sfrac12 & 1 \\
    0 & 1 & 0 & 1 & 0 & \sfrac12
  \end{pmatrix}.
\]
In terms of overall points $A \mathbf{1}$, the winners are the first and third player, tied with $4 \sfrac12$ points and the fourth player is last with $1 \sfrac12$ points. According to Wei's score $A^2 \mathbf{1}$, the first player now gets the highest score with $14 \sfrac14$, whereas the third player only gets $11 \sfrac14$ points. All the other players are tied for last. However, if we continue and compute the score as $A^3 \mathbf{1}$, then the fourth player gets a higher score than players $2$, $5$, and $6$, making him third overall. The fact that the fourth player won against the winner of the tournament gets more and more weight with each step.

The power method is still one of the best methods to compute the largest eigenvalue and a corresponding eigenvector for a large matrix; so it is still used by \emph{Google}, as part of a larger toolbox, to this day.

\section*{The \emph{Google} Matrix}

\begin{figure}
  \includegraphics[width = .66\textwidth]{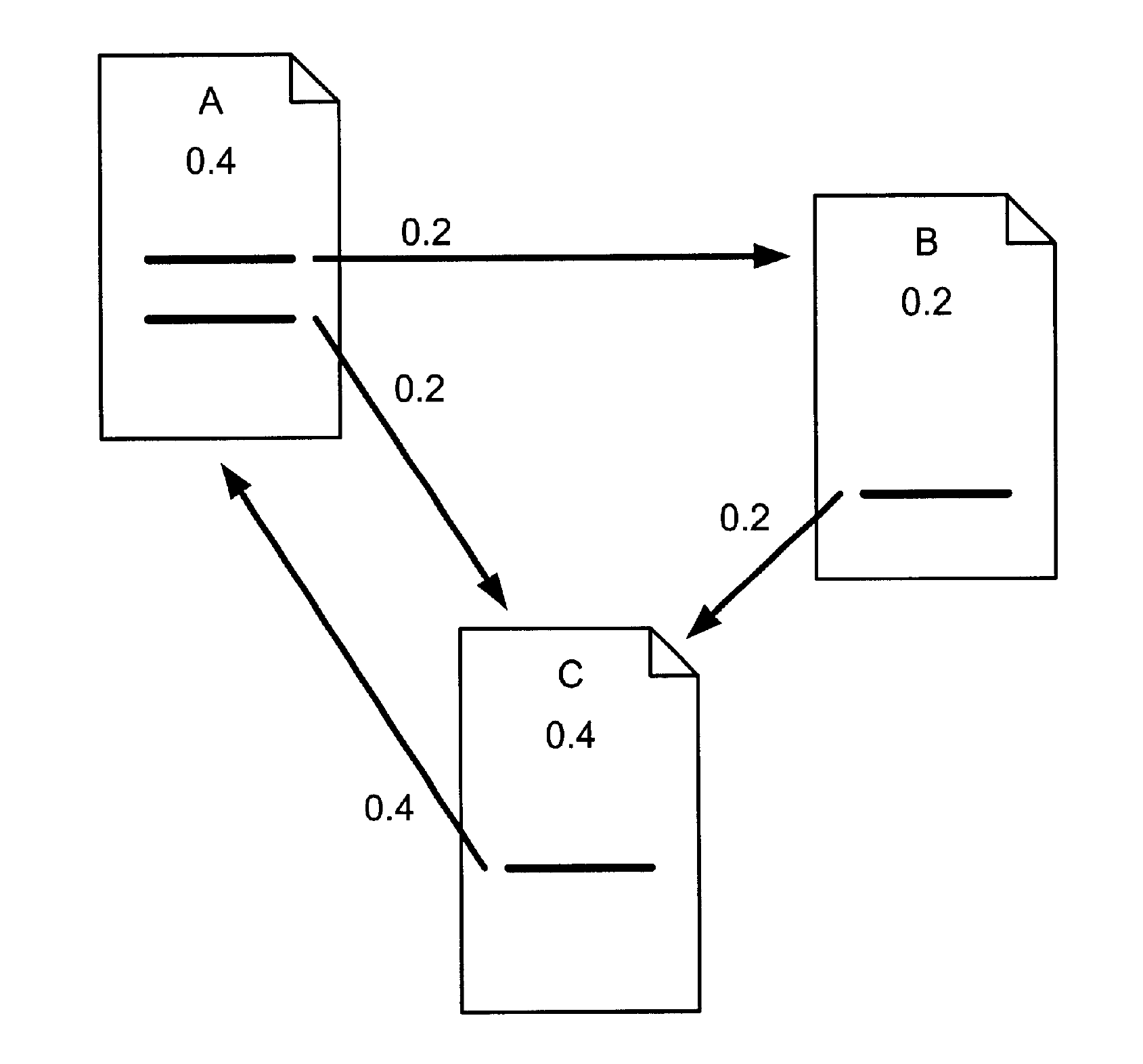}
  \caption{Figure 2 in the \emph{Google} patent US 6,285,999 B1.}
  \label{fig:GooglePatent}
\end{figure}

A commercially amazingly successful development at least partly based on Landau's approach is the original \emph{Google} search algorithm patent, filed in January 1998 by the Stanford PhD student Lawrence “Larry” Page (who started \emph{Google} with his fellow graduate student Sergey Brin), and which was eventually granted to Stanford University in 2001. The central question is what the analogon to a game result 
should be for the internet, formed by a huge number of linked web pages. 
The answer suggested in the patent document is to weigh the relation between websites $i$ and $j$ by the 
fraction of all the hyperlinks on the website $i$ that point to the website $j$. So the inventors of \emph{Google} use the \emph{hyperlink matrix} $H$ with entries 
\[ 
h_{ij} = \left\{
  \begin{array}{cl}
    \frac{t_{ij}}{L_j}, & \mbox{if website $j$ has $t_{ij}$ links to website $i$,} \\
    0, & \mbox{otherwise,}
  \end{array}
  \right.
\]
where $L_j$ denotes the total number of links on the webpage~$j$. 
In the patent application, they give the following example: Suppose we have three website $A$, $B$, and~$C$. Website $A$ has two links, one to $B$ and one to $C$. Website $B$ has only one link to $C$ and Website $C$ has one link to $A$. In that case, the hyperlink matrix $H$ is
\[
  H = \begin{pmatrix}
    0 & 0 & 1 \\
    \sfrac12 & 0 & 0 \\
    \sfrac12 & 1 & 0
  \end{pmatrix}.
\]
The ranks of the websites according to Landau's method are given by solving the eigenvalue problem $H \mathbf{r} = \lambda \mathbf{r}$ for the dominant eigenvalue. In this case, this eigenvalue is $1$ and the corresponding eigenspace gives the ranks $r(A) = 0.4$, $r(B) = 0.2$ and $r(C) = 0.4$ for the three websites, which are the same as the one in the patent application, see \Cref{fig:GooglePatent}.%
\footnote{The matrix $H$ does not have a symmetry any more, which would be analogous to the observation that the off-diagonal entries of $A+A^t$ are $1$ for the Landau matrix, but this is not needed for Landau's method.}

The hyperlink matrix is normalized such that each column sums to $1$. 
So we can interpret the $(i,j)$ entry as the probability to arrive on page $i$ by taking a random link on page $j$. The sum of the entries in row $i$ of this matrix thus reflects the probability to arrive at website $i$ after a random jump from a random website -- this probability is given as the $i$th component of the vector $\sfrac{1}{n} H\mathbf{1}$, whose components are nonnegative and sum to $1$.

This suggests the following new interpretation of the hyperlink matrix~$H$ of Brin and Page: Instead of thinking of the entries of $H$ in terms of pairwise comparisons, let us consider $H$ as the transition matrix of a random walk. Then $H$ describes the probabilities for choices of a surfer browsing a website, who from any webpage would transition to a new one by picking a link on the current website uniformly at random. This point of view in nicely explained in the 1953 paper \cite{Katz} by Leo Katz in the social sciences, which was cited in Page's patent application.
Intuitively, we can think of a website as relatively \enquote{important} if the frequency at which the random walk on the whole internet would visit it is relatively high. 
Thus the \emph{importance} $x_i$ of website $i$ is again, as in Landau's idea, proportional to $H\mathbf{x}$, so that the resulting equation is $H\mathbf{x} =\lambda \mathbf{x}$. However, the internet is not the same as a round robin tournament --- its matrix has a different structure. Not only does it lack the symmetry of a game results matrix, but also other problems could arise. And if they \emph{can} arise, in a large system like the internet they \emph{will}: 

\noindent
\emph{Problem 1}: Not every website has links, so we would have $L_j=0$ in the construction of the hyperlink matrix, and the random walk could get stuck. Brin and Page decided to replace all columns in $H$ where this occurs by $\sfrac{1}{n}\mathbf{1}$, to obtain and then use a modified matrix $S$. The interpretation in terms of the random walk is that if there are no outgoing links on a website, we simply decide uniformly at random where to go next.

\noindent
\emph{Problem 2}: The Perron--Frobenius Theorem for nonnegative matrices has the assumption that the matrix be irreducible. This assumption would be satisfied if one could get from any webpage to any other via a chain of weblinks, that is, if the internet were a connected digraph. In practice, this is not the case.
To get around this problem, Brin and Page proposed a brute force solution, namely to make all entries of the matrix positive (instead of nonnegative), which makes the irreducibility assumption satisfied. Historically, this version of Perron--Frobenius, where all entries of the matrix are positive, was proved first. So the final \emph{Google} matrix is $G = (1-\alpha) S + \alpha \sfrac{1}{n}\,\mathbf{1}_{n\times n}$ for some small but positive value $\alpha$. 
Here, $\mathbf{1}_{n\times n}$ is the $n\times n$ matrix whose entries are all equal to $1$. The \emph{Google} \emph{Page\-Rank} of a website is, just as in Landau's method, the corresponding entry in the essentially unique positive eigenvector corresponding to the dominant, unique, largest eigenvalue a nonnegative square matrix, namely of~$G$.

\section*{Conclusion}

Landau's idea---to calculate the quality of a chess player in a tournament from the strengths of the players he beats or draws with---sounds self-referential, and so it is: There is no way to see or compute the qualities of a single player just from the games he plays. We have to solve the total system, that is, compute the qualities of every player in the tournament.
Thus Landau's innocent-looking first paper leads to the question of self-referentiality (which later in other forms occupied Gödel and Turing), it was later reinterpreted in terms of random walks, it was extended to huge systems that challenged numerical linear algebra (with the power method as a simple but effective solution), and it can still be seen as an ingredient in one of the most impressive economical success stories “powered by mathematics.” Certainly \emph{Google} is powered by many different influences and tools, but it is also “powered by Landau.”

\end{document}